\documentclass[11pt]{amsart}
\usepackage{palatino}
\usepackage{amsfonts}
\usepackage{amssymb}
\usepackage{amscd}

\theoremstyle{plain}

\newtheorem{theorem}{Theorem}[section]
\newtheorem{proposition}[theorem]{Proposition}
\newtheorem{corollary}[theorem]{Corollary}
\newtheorem{lemma}[theorem]{Lemma}
\theoremstyle{definition}

\newtheorem{remark}[theorem]{Remark}

\newtheorem{conjecture}[theorem]{Conjecture}
\newtheorem{conjecture/question}[theorem]{Conjecture/Question}
\newtheorem{question}[theorem]{Question}
\newtheorem{remark/definition}[theorem]{Remark/Definition}
\newtheorem{terminology/notation}[theorem]{Terminology/Notation}

\newcommand{\marginlabel}[1]%
  {\mbox{}\marginpar{\raggedleft\hspace{0pt}\bfseries\sf#1}}

\def\PP{{\textbf P}}
\def\OO{\mathcal{O}}

\def\cB{\mathcal{B}}
\def\cA{\mathcal{A}}

\def\G{\mathcal{G}}

\def\L{\mathcal{L}}

\def\cM{\mathcal{M}}

\def\H{\mathcal{H}}
\def\Pic0{{\rm Pic}^0(X)}

\def\mm{\overline{\mathcal{M}}}
\def\kk{\overline{\mathcal{K}}}

\def\KK{\overline{\mathcal{K}}}

\theoremstyle{remark}

\begin{document}

\title{\bf The global geometry of the moduli space of curves}

\author[G. Farkas]{Gavril Farkas}
\address{Department of Mathematics, University of Texas,
Austin, TX 78712} \email{{\tt gfarkas@math.utexas.edu}}
\thanks{Research  partially supported by the NSF Grant DMS-0450670 and by the Sloan Foundation}

\maketitle

\section{Introduction}

For a complex projective variety $X$,  one way of understanding its
birational geometry is by describing its cones of \emph{ample} and
\emph{effective} divisors
\begin{center}
$\mbox{Ample}(X)\subset \mbox{Eff}(X)\subset N^1(X)_{\mathbb R}.$
\footnote{Throughout this paper we use the formalism of $\mathbb R$-divisors
and we say that a class $e\in N^1(X)_{\mathbb R}$ is
\emph{effective} (resp. \emph{ample}) if $e$ is represented by a
$\mathbb R$-divisor $D$ on $X$ which is effective (resp. ample).}
\end{center}
 The closure in $N^1(X)_{\mathbb R}$ of $\mbox{Ample}(X)$ is the cone
$\mbox{Nef}(X)$ of \emph{numerically effective} divisors, i.e. the set of all
classes $e\in N^1(X)_{\mathbb R}$ such that $C\cdot e\geq 0$ for all curves $C\subset X$. The
interior of the closure $\overline{\mbox{Eff}(X)}$ is the cone of
\emph{big} divisors on $X$. Loosely speaking, one can think of the
nef cone as parametrizing \emph{regular} contractions
\footnote{Recall that a contraction of $X$ is a morphism
$f:X\rightarrow Y$ with connected fibres.} from $X$ to other
projective varieties, whereas the effective cone accounts for
\emph{rational} contractions of $X$. For arbitrary varieties of
dimension $\geq 3$ there is little connection between
$\mbox{Nef}(X)$ and $\mbox{Eff}(X)$ (for surfaces there is Zariski
decomposition which provides a unique way of writing an effective
divisor as a combination of a nef and a "negative" part and this
relates the two cones, see e.g. \cite{L1}). Most questions in higher
dimensional geometry can be phrased in terms of the ample and
effective cones. For instance, a smooth projective variety $X$ is of
general type precisely when $K_X\in
\mbox{int}(\overline{\mbox{Eff}(X)})$.

The question of describing the ample and the effective cone of
$\mm_g$ goes back to Mumford (see e.g. \cite{M1}, \cite{H2}). Moduli
spaces of curves with their inductive structure given by the
boundary stratification are good test cases for many problems coming
from higher dimensional birational geometry. The first major result
result on the global geometry of $\mm_g$ was the celebrated theorem
of Harris, Mumford and Eisenbud that $\mm_g$ is of general type for
all $g\ge 24$ (cf. \cite{HM}, \cite{H1}, \cite{EH3}). This result
disproved a famous conjecture of Severi's who, based on evidence
coming from small genus, predicted that $\cM_g$ is unirational for
all $g$. The space $\mm_g$ being of general type, implies for
instance that the general curve of genus $g\geq 24$ does not appear
in any non-trivial linear system on any non-ruled surface.

The main aim of this paper is to discuss what is currently known
about the ample and the effective cones of $\mm_{g, n}$.
Conjecturally, the ample cone has a very simple description being
dual to the cone spanned by the irreducible components of the locus
in $\mm_{g, n}$ that consists of curves with $3g-4+n$ nodes (cf.
\cite{GKM}). The conjecture has been verified in numerous cases and
it predicts that for large $g$, despite being of general type,
$\mm_g$ behaves from the point of view of Mori theory just like a
Fano variety, in the sense that the cone of curves is polyhedral,
generated by rational curves. In the case of the effective cone the
situation is more complicated. In \cite{FP} we showed that the
Harris-Morrison Slope Conjecture which singled out the Brill-Noether
divisors on $\mm_g$ as those of minimal slope, is false. In this
paper we describe a very general construction of geometric divisors
on $\mm_g$ which provide counterexamples to the Slope Conjecture in
infinitely many genera (see Theorem \ref{div}). Essentially, we
construct an effective divisor of exceptionally small slope on
$\mm_g$ for $g=s(2s+si+i+1)$, where $s\geq 2, i\geq 0$. For $s=1$,
we recover the formula for the class of the Brill-Noether divisor
first computed by Harris and Mumford in \cite{HM}. The divisors
constructed in \cite{EH3}, \cite{FP}, \cite{F2} and \cite{Kh} turn
out to be particular cases of this construction.

 In spite of all the counterexamples, it still seems
reasonable to believe that a "Weak" Slope Conjecture on $\mm_g$
should hold, that is, there should be a universal lower bound
 on the slopes of all effective divisors on $\mm_g$ which is independent of
$g$. This fact would highlight a fundamental difference between
$\mathcal{M}_g$ and $\mathcal{A}_g$ and would provide a modern
solution to the \emph{Schottky problem} (see Subsection 2.2 for more
details). In Section 3 we announce a proof that $\mm_{22}$ is of
general type and we describe the Kodaira type of the moduli spaces
$\mm_{g, n}$ of $n$-pointed stable curves.

\section{Divisors on $\mm_{g,n}$}

For non-negative integers $g$ and $n$ such that $2g-2+n>0$ we denote by $\mm_{g, n}$ the
moduli stack of $n$-pointed stable curves of genus $g$. The stack $\mm_{g, n}$ has a
stratification given by \emph{topological type}, the codimension $k$
strata being the components of the closure in $\mm_{g, n}$ of the
locus of curves with $k$ nodes. The $1$-dimensional strata are also
called \emph{$F$-curves} and it is easy to see that each $F$-curve
is isomorphic to either  $\mm_{0, 4}$ or to $\mm_{1,1}$. It is
straightforward to list all $F$-curves on a given $\mm_{g, n}$. For
instance, $F$-curves on $\mm_{0, n}$ are in $1:1$ correspondence
with partitions $(n_1, n_2, n_3, n_4)$ of $n$, the corresponding
$F$-curve being the image of the map $\nu:\mm_{0, 4}\rightarrow
\mm_{0, n}$ which takes a rational $4$-pointed curve $[R, x_1, x_2,
x_3, x_4]$ to a rational $n$-pointed curve obtained by attaching  a
fixed rational $(n_i+1)$-pointed curve at the point $x_i$.

The codimension $1$ strata in the topological stratification are the boundary  divisors on $\mm_{g,n}$
which are indexed as follows: For $0\leq i\leq g$ and $S\subset \{1,\ldots,n\}$, we denote by
 $\delta_{i:S}$ the class of the closure of the locus of nodal curves
 $C_1\cup C_2$, where $C_1$ is a smooth curve of genus $i$, $C_2$ is a
 smooth curve
 of genus $g-i$ and such that the marked points sitting on
 $C_1$ are precisely those labeled by $S$.
We also have the class $\delta_{0}$ corresponding to irreducible
pointed curves with a single node. Apart from boundary divisor
classes, we also introduce the tautological classes
$\{\psi_i=c_1(\mathbb L_i)\}_{i=1}^n$  corresponding to the $n$
marked points. Here $\mathbb L_i$ is the line bundle over the moduli
stack with fibre $\mathbb L_i[C, x_1, \ldots,
x_n]:=T_{x_i}(C)^{\vee}$ over each point $[C, x_1, \ldots, x_n]\in
\mm_{g, n}$. Finally, we have the Hodge class defined as follows: If
$\pi:\mm_{g, 1}\rightarrow \mm_g$ is the universal curve, we set
$\lambda:=c_1(\mathbb E)\in \mbox{Pic}(\mm_g)$, where $\mathbb
E:=\pi_*(\omega_{\pi})$ is the rank $g$ Hodge bundle on $\mm_g$.
 A result of Harer and Arbarello-Cornalba (cf. \cite{AC}) says
 that $\lambda, \psi_1,\ldots, \psi_n$ together with the boundary classes  $\delta_{0}$ and $\delta_{i:S}$ generate
 $\mbox{Pic}(\overline{\mathcal{M}}_{g,n})$. This result has been extended to arbitrary characteristic by Moriwaki (cf. \cite{Mo2}). When $g\geq 3$ these classes
 form a basis of
  $\mbox{Pic}(\mm_{g, n})$.

\subsection{The ample cone of $\mm_{g, n}$}\hfill

In this section we describe the ample cone of $\mm_{g,n}$.
Historically speaking, the study of ample divisors on $\mm_g$ began when
Cornalba and Harris  proved that  the $\mathbb Q$-
class $a\lambda-\delta_0 -\cdots - \delta_{[g/2]}$ is ample on
$\mm_g$ if and only if $a>11$ (cf. \cite{CH}). Later, Faber
completely determined $\mbox{Ample}(\mm_3)$: a class $D\equiv
a\lambda -b_0\delta_0-b_1\delta_1 \in \mbox{Pic}(\mm_3)$ is nef if
and only if
$$ 2b_0-b_1\geq 0, \mbox{ }\mbox{ }b_1\geq 0 \mbox{ and } \mbox{ }a-12b_0+b_1\geq
0\mbox{ } \mbox{ } \mbox{(cf. \cite{Fa})}.$$ He pointed out that the numbers
appearing in the left hand side of these inequalities are
intersection numbers of $D$ with certain $F$-curves in $\mm_3$ thus
raising for the first time the possibility that the $F$-curves might
generate the Mori cone of curves $NE_1(\mm_{g, n})$. The
breakthrough in this problem came when Gibney, Keel and Morrison
proved that strikingly, $NE_1(\mm_{g, n})$ is the sum of the cone
generated by $F$-curves and  the cone $NE_1(\mm_{0, g+n})$. In this way,
computing the nef cone of $\mm_{g,n}$ for any $g>0$ always boils
down to a problem in genus $0$!

\begin{theorem}\label{bridge}([GKM])
If $j:\mm_{0, g+n}\rightarrow \mm_{g,n}$ is the \lq \lq flag map"
given by attaching fixed elliptic tails to the first $g$ marked
points of every $(g+n)$-pointed stable rational curve, then a divisor $D$
on $\mm_{g,n}$ is nef if and only if $j^*(D)$ is nef on $\mm_{0,
g+n}$ and $D\cdot C\geq 0$ for every $F$-curve $C$ in $\mm_{g,n}$.
\end{theorem}

This reduction to genus $0$ then makes the following conjecture very plausible:
\begin{conjecture}([GKM])\label{GKM}
The Mori cone $NE_1(\mm_{g,n})$ is generated by $F$-curves. A
divisor $D$ on $\mm_{g, n}$ is ample if and only if $D\cdot C>0$ for
every $F$-curve $C$ in $\mm_{g,n}$.
\end{conjecture}

The conjecture reflects the expectation that the extremal rays of
$\mm_{g, n}$ should have modular meaning. Since $F$-curves can be
easily listed, this provides an explicit (conjectural) description
of the ample cone. For instance, on $\mm_g$, the conjecture predicts
that a divisor $D\equiv
a\lambda-b_0\delta_0-\cdots-b_{[g/2]}\delta_{[g/2]}\in
\mbox{Pic}(\mm_g)$ is ample if and only if the following
inequalities are satisfied:
$$a-12b_0+b_1> 0, \mbox{ }\mbox{ } 2b_0> b_i> 0 \mbox{ for all }i\geq 1,
\mbox{ } \mbox{ }b_i+b_j> b_{i+j} \mbox{ for all }i, j\geq 1 \mbox{ with }i+j\leq g-1$$ and
$$b_i+b_j+b_k+b_l> b_{i+j}+b_{i+k}+b_{i+l} \mbox{ for all }i, j,
k, l\geq 1 \mbox{ with }i+j+k+l=g.$$ Here we have the usual
convention $b_i=b_{g-i}$. Conjecture \ref{GKM} has been checked on
$\mm_g$ for all $g\leq 24$ in \cite{KMcK}, \cite{FG} and \cite{G}.
In fact, Gibney has reduced the conjecture on a given $\mm_g$ to an
entirely combinatorial question which can be checked by computer.
Recently, Coskun, Harris and Starr have reduced the calculation of
the ample cone of the moduli space of stable maps
$\mm_{0,n}(\PP^r,d)$ to Conjecture \ref{GKM} for $\mm_{0, n+d}$ (cf.
\cite{CHS}). In \cite{GKM} it is also pointed out that Conjecture
\ref{GKM} would be implied by an older conjecture of Fulton
motivated by the fact that $\mm_{0,n}$ has many of the geometric features of a
toric variety (without being a toric variety, of course):

\begin{conjecture}
Any divisor class $D$ on $\mm_{0, n}$ satisfying $C\cdot D\geq 0$
for all $F$-curves $C$ in $\mm_{0, n}$ can be expressed as an
effective combination of boundary classes.
\end{conjecture}

Fulton's conjecture is true on $\mm_{0,n}$ for $n\leq 6$ (cf.
\cite{FG}). Note that it is not true that every effective divisor on
$\mm_{0, n}$ is equivalent to an effective combination of boundary
divisors (cf. \cite{Ve}): If $\xi:\mm_{0, 6}\rightarrow \mm_3$
denotes the map which identifies three pairs of marked points on a
genus $0$ curve, then the pull-back of the hyperelliptic locus
$\xi^*(\mm_{3, 2}^1)$ is an effective divisor on $\mm_{0, 6}$ for
which there exists an explicit curve $R\subset \mm_{0, 6}$ not
contained in the boundary of $\mm_{0, 6}$ such that $R\cdot
\xi^*(\mm_{3, 2}^1)<0$. Thus $\xi^*(\mm_{0, 6})$ is not an effective
combination of boundary divisors.

\begin{remark} In low genus one can show that $\mbox{Ample}(\mm_g)$ is \lq
 \lq tiny" inside the much bigger cone $\mbox{Eff}(\mm_g)$ which
 underlies the fact that regular contractions of $\mm_g$ do not capture the rich birational geometry of $\mm_g$ (For instance, the only divisorial contraction
 of $\mm_{g, n}$ with relative Picard number $1$ is the blow-down of the elliptic tails, see
 \cite{GKM}). The difference between the two
 cones can be vividly illustrated on $\mm_3$: we have seen that
 $\mbox{Nef}(\mm_3)$ is generated by the classes $\lambda,
 12\lambda-\delta _0$ and $10\lambda-\delta_0-2\delta_1$ (cf. \cite{Fa}), whereas it is easy to show that $\mbox{Eff}(\mm_3)$ is much larger, being spanned by  $\delta_0, \delta_1$ and the class of the hyperelliptic locus
 $h=9\lambda-\delta_0-3\delta_1$.
 \end{remark}

 Theorem \ref{bridge} has a number of important applications to the
 study of regular morphisms from $\mm_{g, n}$ to other projective
 varieties. For instance it is known that for $g\geq 2$, $\mm_g$ has
 no non-trivial fibrations (that is, morphisms with connected fibres
 to lower dimensional varieties). Any fibration of
 $\mm_{g, n}$ must factor through one of the forgetful maps $\mm_{g,
 n}\rightarrow \mm_{g, i}$ for some $i<n$ (see \cite{GKM}, Corollary
 0.10).  If  $f:\mm_{g,
 n}\rightarrow X$ is  a birational morphism to a projective variety, it is known that the
 exceptional locus $\mbox{Exc}(f)$ is contained in the boundary of
 $\mm_{g, n}$. In particular such a projective variety $X$ is a new
 compactification of $\cM_g$. (The use of such a result is
 limited however by the fact that there are very few known examples
 of regular morphism from $\mm_{g, n}$). Theorem \ref{bridge} can
 be directly applied to show that many types of divisors $D$ on
 $\mm_{g, n}$ which non-negatively meet all $F$-curves are actually
 nef. For instance one has the following result (cf. \cite{GKM},
 Proposition 6.1):

 \begin{theorem}\label{gkm2}
 If $D\equiv a \lambda -\sum_{i=0}^{[g/2]} b_i \delta_i$ is a
 divisor on $\mm_g$ such that $b_i\geq b_0$ for all $1\leq i\leq [g/2]$ and
 $C\cdot D\geq 0$ for any $F$-curve $C$, then $D$ is nef.
 \end{theorem}

\begin{remark} Since any regular morphism $f:\mm_{g, n}\rightarrow
X$ to a projective variety is responsible for a semi-ample line
bundle $L:=f^*(\OO_X(1))$ rather than a nef one, it is a very
interesting question to try to characterize semi-ample line bundles
on $\mm_{g, n}$. Surprisingly, this question is easier to handle in
positive characteristic due to the following result of Keel (cf.
\cite{K}, Theorem 0.2): If $L$ is  a nef line bundle on a projective
variety $X$ over a field of positive characteristic, then $L$ is
semi-ample if and only if the restriction of $L$ to its
\emph{exceptional locus}  $\mbox{Exc(L)}$ is semi-ample. Recall
that if $L$ is a nef line bundle on $X$, then $$
\mbox{Exc}(L):=\bigcup \{Z\subset X: Z \mbox{ is an irreducible
subvariety with
 }L^{\mbox{dim}(Z)}\cdot Z=0\}.$$
 An easy application of Keel's Theorem is that the
tautological class $\psi\in \mbox{Pic}(\mm_{g, 1})$ is semi-ample in
positive characteristic but fails to be so in characteristic $0$
(see \cite{K}, Corollary 3.1). No example of a nef line bundle on
$\mm_g$ which fails to be semi-ample is known although it is
expected there are many such examples.
\end{remark}

 \begin{remark} It makes sense of course to ask what is the nef cone
 of the moduli space of abelian varieties. Shepherd-Barron computed the nef
 cone of the first Voronoi compactification $\mathcal{A}_g^{I}$ of $\mathcal{A}_g$ (cf.
 \cite{SB}). Precisely, $NE_1(\mathcal{A}_g^{I})$ is generated by two curve classes $C_1$ and
 $C_2$,
 where $C_1$ is any exceptional curve in the contraction of $\mathcal{A}_g^{I}$ to the Satake compactification
 of $\mathcal{A}_g$, while $C_2=\{[X\times E]\}_{[E]\in \mathcal{A}_1}$, where $[X]\in \mathcal{A}_{g-1}$
 is a fixed ppav of dimension $g-1$ and $E$ is a moving elliptic curve. Hulek and Sankaran have
 determined the nef cone of the second Voronoi
 compactification $\mathcal{A}_4^{II}$ of $\mathcal{A}_4$ (cf. \cite{HS}).
 \end{remark}

\noindent {\bf{Towards the canonical model of $\mm_{g, n}$}} \vskip
5pt \noindent
 A somewhat related question concerns the canonical
model of $\mm_g$. Since the variety $\mm_g$ is of general type for
large $g$, a result from \cite{BCHM} implies the finite generation
of the canonical ring $R(\mm_g):=\oplus_{n\geq 0}H^0(\mm_g, n
K_{\mm_g})$ and the existence of a canonical model of the moduli
space $\mm_g^{can}:=\mbox{Proj}\bigl(R(\mm_g)\bigr)$. It is natural
to ask for a modular interpretation of the canonical model. Very
interesting ongoing work of Hassett and Hyeon provides the first
steps towards understanding $\mm_g^{can}$ (see \cite{Ha1},
\cite{HH}, but also \cite{HL} where the Minimal Model Program for
$\mm_3$ is completed). Precisely, if
$$\delta:=\delta_0+\cdots+\delta_{[g/2]}\in \mathrm{Pic}(\mm_g)$$ denotes the total boundary
of $\mm_g$ and $K_{\overline{M}_g}=13\lambda-2\delta$ is the
canonical class of the moduli stack, for each rational number $0\leq
\alpha \leq 1$ we introduce the \emph{log canonical model}
$$\mm_g^{can}(\alpha):=\mbox{Proj}\Bigl(\oplus_{n\geq 0}
H^0\bigl(\mm_g, n(K_{\overline{M}_g}+\alpha \delta)\bigr)\Bigr).$$
Then $\mm_g^{can}(\alpha)=\mm_g$ for $9/11\leq \alpha \leq 1$
because of the already mentioned result of Cornalba and Harris
\cite{CH}, whereas $\lim_{\alpha\rightarrow 0}
\mm_g^{can}(\alpha)=\mm_g^{can}$. The first interesting question is
what happens to $\mm_g^{can}(\alpha)$ when $\alpha=9/11$ since in
this case there there exists a curve $R\subset \mm_g$ such that
$(K_{\overline{M}_g}+\frac{9}{11} \delta)\cdot R=0$ (precisely, $R$
corresponds to a pencil of plane cubics with a section which is
attached to a fixed pointed curve of genus $g-1$). It turns out that
for $7/10< \alpha \leq 9/11$, the moduli space $\mm_g^{can}(\alpha)$
exists and it is identified with the space $\mm_g^{\mathrm{ps}}$ of
\emph{pseudo-stable} curves in which cusps are allowed but elliptic
tails are ruled out. The morphism
$\mm_g\stackrel{|11\lambda-\delta|}\longrightarrow
\mm_g^{\mathrm{ps}}$ is a divisorial contraction of the boundary
divisor $\Delta_1$. The next (substantially more involved) step is
to understand what happens when $\alpha=7/10$. It turns out that
$\mm_g^{can}(7/10)$ exists as the quotient by $SL_{3g-3}$ of the
Chow variety of bicanonical curves $C\subset \PP^{3g-4}$, whereas
the model $\mm_g^{can}(7/10-\epsilon)$ for $0\leq \epsilon <<1$
exists and it is obtained from $\mm_g^{can}(7/10+\epsilon)$ by an
explicit flip and it parameterizes curves with nodes, cusps and
tacnodes as singularities. As $\alpha \rightarrow 0$, one expects of
course worse and worse singularities like higher-order tacnodes.
\subsection{The effective cone of $\mm_{g}$}\hfill
Following Harris and Morrison \cite{HMo}, we define the slope
function on the effective cone $s:\mbox{Eff}(\mm_g) \rightarrow
\mathbb R\cup \{\infty\}$ by the formula
$$s(D):= \mbox{inf }\{\frac{a}{b}:a,b>0 \mbox{ such that }
a\lambda-b\delta-D\equiv \sum_{i=0}^{[g/2]} c_i\delta_i,\mbox{ where
}c_i\geq 0\}.$$ From the definition it follows that $s(D)=\infty$
unless $D\equiv a\lambda-\sum_{i=0}^{[g/2]} b_i\delta_i$ with
$a,b_i\geq 0$ for all $i$. Moreover, it is well-known that
$s(D)<\infty$ for any $D$ which is the closure of an effective
divisor on $\mathcal{M}_g$. In this case one has that
$s(D)=a/{\mbox{min}}_{i=0}^{[g/2]}b_i.$ We denote by $s(\mm_g)$ the
\emph{slope of the moduli space} $\mm_g$, defined as
$s(\mm_g):=\mbox{inf }\{s(D):D\in \mbox{Eff}(\mm_g)\}$.

\begin{conjecture} (Harris,  Morrison, \cite{HMo})\label{slope}
We have that $s(D)\geq6+12/(g+1)$ for all $D\in \mbox{Eff}(\mm_g)$, with equality if $g+1$ is
 composite and $D$ is a combination of Brill-Noether divisors.
\end{conjecture}

Let us recall the definition of the classical Brill-Noether
divisors. We fix a genus $g\geq 3$ such that there exist $r, d\geq
1$ with $\rho(g, r, d):=g-(r+1)(g-d+r)=-1$ (in particular, $g+1$ has
to be composite). We define the geometric subvariety of $\cM_g$
$\cM_{g, d}^r:=\{[C]\in \cM_g: C \mbox{ has a linear series of type
}\mathfrak g^r_d\}.$ Since $\rho(g, r, d)$ is the expected dimension
of the determinantal variety $W^r_d(C)$ of $\mathfrak g^r_d$'s on a
fixed curve $C$ of genus $g$ (see \cite{ACGH}), one would naively
expect that $\cM_{g, d}^r$ is a divisor on $\cM_g$. In fact we have
a stronger result due to Eisenbud and Harris (cf. \cite{EH2},
\cite{EH3}):

\begin{theorem}
The locus $\cM_{g, d}^r$ is an irreducible divisor on $\cM_g$
whenever $\rho(g, r, d)=-1$. Moreover, the class of the
compactification in $\mm_g$ of the Brill-Noether divisor is given by
the formula
$$\mm_{g, d}^r\equiv c_{g, d, r}\Bigl((g+3)\lambda-\frac{g+1}{6}
\ \delta_0-\sum_{i=1}^{[g/2]} i(g-i)\delta_i \Bigr),$$ where $c_{g,
d,r}$ is an explicitly given constant.
\end{theorem}

 Thus $s(\mm_{g, d}^r)=6+12/(g+1)$ and then the Slope
 Conjecture singles out the Brill-Noether divisors on $\mm_g$ as those
 having minimal slope. Apart from the evidence coming from low genus, the conjecture was mainly
 based on  the large number of calculations of classes of geometric divisors on $\mm_g$
 (see e.g. \cite{EH3}, \cite{H1}).
\begin{remark}
If $D\equiv a \lambda -b_0\delta_0-\cdots -b_{[g/2]} \delta_{[g/2]}$
is a nef divisor on $\mm_g$, then combining the inequalities
$a-12b_0+b_1\geq 0$ and $2b_0-b_1\geq 0$ obtained by intersecting
$D$ with two $F$-curves, we obtain that $s(D)\geq a/b_0\geq 10$. On
the other hand, the class $10\lambda-\delta_0-2\delta_1$ is nef on
$\mm_g$ (cf. \cite{GKM}), hence $\mbox{liminf}\{ s(D): D\in
\mbox{Nef}(\mm_g)\}=10$. This once again illustrates that nef
divisors contribute little to the birational geometry of $\mm_g$.
\end{remark}
As explained in the original paper \cite{HMo}, the Slope Conjecture
is intimately related to the problem of determining the Kodaira
dimension of $\cM_g$. Recall first the computation of the canonical
class of $\mm_g$:

\begin{theorem}\label{canonical}(\mbox{Harris-Mumford})
$$K_{\mm_g}\equiv
13\lambda-2\delta_0-3\delta_1-2\delta_2-\cdots-2\delta_{[g/2]}.$$
\end{theorem}
\begin{proof} If $\pi:\mm_{g, 1}\rightarrow \mm_g$ is the universal
curve, then by Kodaira-Spencer theory we have that
$\Omega_{\mm_g}^1=\pi_*(\omega_{\pi}\otimes \Omega_{\pi})$, where
$\Omega_{\pi}$ is the cotangent sheaf while $\omega_{\pi}$ is the
dualizing sheaf. Then apply Grothendieck-Riemann-Roch to the
universal curve, to obtain that the canonical class of the
\emph{moduli stack} is
$K_{\overline{M}_g}=13\lambda-2(\delta+\ldots+\delta_{[g/2]})$. To
obtain the formula for the canonical class $K_{\mm_g}$ of the \emph{
coarse moduli space} we use that the natural map from the stack to
the coarse moduli space is simply branched along the boundary
divisor $\Delta_1$.
\end{proof}

Since the Hodge class $\lambda$ is big and nef, it follows that
$\mm_g$ is of general type whenever $s(\mm_g)<s(K_{\mm_g})=13/2$.
Since $s(\mm_{g, d}^r)=6+12/(g+1) <13/2 \Longleftrightarrow g\geq
24$, we obtain the main result from \cite{HM} and \cite{EH3}, namely
that $\mm_g$ is of general type for $g\geq 24$ (Strictly speaking,
this argument works only for those $g$ for which $g+1$ is composite.
In the remaining cases, when there are no Brill-Noether divisors on
$\mm_g$, one has to use the locus where the classical Petri Theorem
fails, see \cite{EH3}). The Slope Conjecture would immediately imply
the following statement:

\begin{conjecture}\label{gentype}
The Kodaira dimension of $\mm_g$ is $-\infty$ for all $g<23$.
\end{conjecture}

\noindent {\bf{The unirationality of $\mm_{14}$}}

Severi proved that $\mm_g$ is unirational for $g\leq 10$. The cases
$g=11, 12, 13$ were settled by Sernesi and then Chang and Ran (cf.
\cite{Se}, \cite{CR1}). Moreover, it is known that $\mm_{15}$ is
rationally connected (cf. \cite{BV}) and that
$\kappa(\mm_{16})=-\infty$ (cf. \cite{CR2}). Optimal bounds for
rationality of $\mm_{g, n}$ when $g\leq 6$ are provided in
\cite{CF}. Verra has recently settled the long standing case of
$\mm_{14}$ proving the following theorem (cf. \cite{Ver}):

\begin{theorem}\label{verra}
The moduli space $\cM_{14}$ is unirational.
\end{theorem}
\noindent{\emph{Sketch of proof}.}
 We denote by $H_{d, g, r}$ the Hilbert scheme of curves $C\subset
\PP^r$ with $g(C)=g, \mbox{deg}(C)=d$. The key observation is that
if $[C]\in \cM_{14}$ is suitably general, then $\mbox{dim}
W^6_{18}(C)=0$ and if $C\stackrel{|L|}\hookrightarrow \PP^6$ is the
embedding given by any $L\in W^6_{18}(C)$, then $$\mbox{dim} \
\mbox{Ker}\{\mbox{Sym}^2 H^0(C, L)\rightarrow H^0(C, L^{\otimes
2})\}=5,$$ that is, $C$ lies precisely on $5$ quadrics $Q_1, \ldots,
Q_5$. Writing $\cap_{i=1}^5 Q_i=C\cup R$, one sees that the residual
curve $R$ is smooth with $\mbox{deg}(R)=14$ and $g(R)=8$.
Significantly, $h^1(\OO_R(1))=0$ and the Hilbert scheme $H_{14, 8,
6}$ is a lot easier to study than the scheme $H_{18, 14, 6}$ we
started with and which parametrizes curves $C\subset \PP^6$ with
$h^1(\OO_C(1))=2$. Using Mukai's result that a generic canonical
curve of genus $8$ is a linear section of $G(2, 6)\subset \PP^{14}$,
one proves that $H_{14, 8, 6}$ is unirational. If
$\mathcal{G}_5\rightarrow H_{14, 8, 6}$ denotes the Grassmann bundle
consisting of pairs $[R, V]$ with $[R]\in H_{14, 8, 6}$ and $V\in
G\bigl(5, H^0(I_R(2))\bigr)$, then $\mathcal{G}_5$ is unirational
(because $H_{14, 8, 6}$ is so),  and there exists a dominant
rational map $\mathcal{G}_5 --> H_{18, 14, 6}$ which sends $[R, V]$
to $[C]$, where $C\cup R$ is the scheme in $\PP^6$ defined by $V$.
By standard Brill-Noether theory, the forgetful morphism $H_{18, 14,
6}-->\cM_{14}$ is dominant, hence the composition
$\mathcal{G}_5-->\cM_{14}$ is dominant as well. This shows that
$\cM_{14}$ is unirational. \hfill $\Box$

The Slope Conjecture is also connected to the Schottky problem of
describing geometrically $\cM_g$ in the Torelli embedding
$t:\cM_g\rightarrow \cA_g$ given by $[C]\mapsto [\mbox{Jac}(C),
\Theta_C]$. The map $t$ can be extended to a rational map $t:\mm_g
-->\mathcal{A}_g^{part}$ well-defined at least in codimension $1$,
where $\mathcal{A}_g^{part}$ is Mumford's partial compactification
of rank $1$ degenerations obtained by blowing-up the open subvariety
$\mathcal{A}_g \cup \mathcal{A}_{g-1}$ inside the Satake
compactification of $\mathcal{A}_g$ (cf. \cite{M2}).  One has that
$\mbox{Pic}(\mathcal{A}_g^{part})\otimes \mathbb Q=\mathbb Q\cdot
\lambda\oplus \mathbb Q\cdot \delta$, where $\lambda:=c_1(\mathbb
E)$ is the Hodge class corresponding to modular forms of weight one
and $\delta=[\mathcal{A}_g^{part}-\mathcal{A}_g]$ is the class of
the irreducible boundary divisor. Note that $t^*(\lambda)=\lambda\in
\mbox{Pic}(\mm_g)$ while $t^*(\delta)=\delta_0$. The
quasi-projective variety $\mathcal{A}_g^{part}$ is as good as any
projective model of $\mathcal{A}_g$ when it comes to codimension $1$
problems like determining the Kodaira dimension of of
$\mathcal{A}_g^{part}$ or describing
$\mbox{Eff}(\mathcal{A}_g^{part})$. In particular, one can define
the \emph{slope} of $\mathcal{A}_g$ as being
$$s(\mathcal{A}_g)=s(\mathcal{A}_g^{part}):=\mbox{inf}\{s(D)=\frac{a}{b}: D\equiv
a\lambda-b\ \delta \in \mbox{Eff}(\mathcal{A}_g^{part})\}.$$

\begin{theorem}(\mbox{Tai, [T]})
\
We have that $\mbox{lim}_{g\rightarrow \infty} s(\mathcal{A}_g)=0$.
\end{theorem}

If we combine Tai's estimate with the Slope Conjecture, it would
follow that any Siegel modular form of slope $<6+12/(g+1)$ would
automatically vanish on $\mathcal{M}_g$ thus providing a Schottky
relation. Note that any weaker estimate of the form $s(\mm_g)\geq
\epsilon$ for $g$ large, where $\epsilon>0$ is  a constant
independent on $g$, would suffice to obtain the same conclusion. It
is then very tempting to ask whether the modular forms of slope
$\geq \epsilon$ cut out precisely $\mathcal{M}_g$. A positive answer
to this question would represent a completely novel solution to the
Schottky problem.

Unfortunately, the Slope Conjecture (at least in its original form),
turns out to be false. The first counterexample was constructed in
\cite{FP} and we start by giving a geometric reinterpretation to the
Slope Conjecture which will turn out to be crucial in constructing
counterexamples:

\begin{proposition}\label{nef}
Let $D$ be an effective divisor on $\mm_g$. If $s(D)<6+12/(g+1)$,
then $D$ contains the closure of the locus $\mathcal{K}_g:=\{[C]\in
\cM_g: C\mbox{ sits on a }K3 \mbox{ surface}\}$.
\end{proposition}
\begin{proof}
Clearly, we may assume that $D$ is the closure of an effective
divisor on $\cM_g$. We consider a Lefschetz pencil of curves of
genus $g$ lying on a general $K3$ surface of degree $2g-2$ in
$\PP^g$. This gives rise to a curve $B$ in the moduli space
$\overline{\mathcal{M}}_g$. Since $\mathcal{K}_g$ is the image of a
$\PP^g$-bundle over the irreducible moduli space of polarized $K3$
surfaces of degree $2g-2$, the pencils $B$ fill up the entire locus
$\mathcal{K}_g$. We have that $\lambda\cdot B=g+1$, $ \delta_0 \cdot
B=6g+18$ and $\delta_i\cdot B=0$ for $i\geq 1$. The first two
intersection numbers are computed using the classical formula for
the number of singular fibres in a pencil (see e.g. \cite{GH}, p.
508) while the last assertion is obvious since a Lefschetz pencil
contains no reducible curves. We can write that
 $ \delta \cdot B / \lambda \cdot B=6+12/(g+1)>s(D)$, which implies that $D \cdot B<0$ hence $B\subset D$.
  By varying $B$ and $S$ we obtain that $\KK_g\subset D$.
\end{proof}

\begin{remark}
 Proposition \ref{nef} shows  that the Slope Conjecture would be implied by the curve
 $B\subset \mm_g$ being nef. It is also proved in \cite{FP} that
 if $D\equiv a\lambda-b_0\delta_0-\cdots -b_{[g/2]}\delta_{[g/2]}\in
 \mbox{Eff}(\mm_g)$ is such that $a/b_0\leq 71/10$, then $b_i\geq b_0$ for all $1\leq i\leq 11$.
 At least for $g\leq 23$, the statement of the Slope Conjecture is
 thus equivalent to $B$ being a nef curve.
\end{remark}

\begin{remark}
The pencils $B$ fill up $\mm_g$ for $g\leq 11, g\neq 10$ (cf.
\cite{Mu}), hence Proposition \ref{nef} gives a short proof of the
Slope Conjecture on $\mm_g$ for these values.  For those $g$ such
that $\mathcal{K}_g\subsetneq \cM_g$, Proposition \ref{nef} suggests
how to search for counterexamples to the Slope Conjecture: one has
to come up with divisorial geometric properties which are a
relaxation of the condition that a curve lie on a $K3$ surface. The
first case where one can test the Slope Conjecture is on $\mm_{10}$,
where contrary to the naive dimension count, $\mathcal{K}_{10}$ is a
divisor: The moduli space of polarized $K3$ surfaces of genus $g$
depends on $19$ parameters, hence the expected dimension of
$\mathcal{K}_g$ is $\mbox{min}(19+g, 3g-3)$ which would suggest that
any $[C]\in \cM_{10}$ lies on a $K3$ surface.  However, Mukai has
proved that $K3$ surfaces of genus $10$ appear as codimension $3$
linear sections of a certain rational homogeneous variety
$X_5\subset \PP^{13}$ corresponding to the Lie group $G_2$ (cf.
\cite{Mu}). Therefore, if $[C]\in \cM_{10}$ lies on a $K3$ surface,
then $C$ lies on $\infty^3$ $K3$ surfaces and $\mathcal{K}_{10}$ is
a divisor on $\mm_{10}$. The Slope Conjecture holds on $\mm_{10}$ if
and only if it holds for $\KK _{10}$.
\end{remark}

\begin{theorem}([\mbox{FP}])\label{class} The divisor $\KK_{10}$ provides a
counterexample to the Slope Conjecture. Its class is given by the
formula $\KK_{10}\equiv
7\lambda-\delta_0-5\delta_1-9\delta_2-12\delta_3-14\delta_4-15\delta_5$,
hence $s(\KK_{10})=7$.
\end{theorem}

The proof of this theorem does not use the original definition of
$\KK_{10}$. Instead, we show that $\KK_{10}$ has a number of other
interpretations, in particular, we can geometrically characterize
the points from $\kk_{10}$ in ways that make no reference to $K3$
surfaces and use these descriptions to compute the class of
$\KK_{10}$.

\begin{theorem}([\mbox{FP}])\label{equiv}
The divisor $\mathcal{K}_{10}$ has two other incarnations as a geometric
subvariety of $\mathcal{M}_{10}$:
\begin{enumerate}
\item The locus of curves $[C]\in \cM_{10}$ carrying  a semistable rank two vector
 bundle $E$ with $\wedge^2(E)=K_C$ and $h^0(C, E)\geq 7$.
\item The locus of curves $[C]\in \cM_{10}$ for which there exists
$L\in W^4_{12}(C)$ such that the multiplication map
$\rm{Sym}$$^2H^0(L)\rightarrow H^0(L^{\otimes 2})$ is not an
isomorphism.
\end{enumerate}
\end{theorem}

\noindent \emph{Proof of Theorem \ref{class}.} We use the second
description from Theorem \ref{equiv} and studying degenerations of
multiplication maps, we can explicitly describe the pull-backs
$j_i^*(\KK_{10})$ where $j_i^*:\mm_{i, 1}\rightarrow \mm_{10}$ is
the map obtained by attaching a fixed tail of genus $10-i$ at the
marked point of each genus $i$ curve. It turns out that these
pull-backs are sums of \lq \lq pointed" Brill-Noether divisors on
$\mm_{i, 1}$. Since these classes have been computed in \cite{EH2},
we get enough relations in the coefficients of $[\KK_{10}]$ that
enable us to fully determine its class.\hfill $\Box$

\begin{remark} Note that if $B\subset \mm_{10}$ is the pencil appearing in  Proposition
\ref{nef}, then $\kk_{10}\cdot B=-1$. The Slope Conjecture fails on $\mm_{10}$ precisely because of the
failure of $B$ to be a nef curve.
\end{remark}

\begin{remark} We note that a general curve $[C]\in \cM_{10}$ possesses
finitely many (precisely $42$) linear series $\mathfrak
g^4_{12}=K_C(-\mathfrak g^1_6)$, and these $\mathfrak g^1_6$'s are
the pencils of minimal degree on $C$. If $C$ lies on a $K3$ surface
$S$, the exceptional rank $2$ vector bundle $E$ which appears in
Theorem \ref{equiv} is a Lazarsfeld-Mukai bundle obtained as a
restriction to $C$ of a rank $2$ bundle on $S$ which is the
elementary transformation along $C$ given by the evaluation map
$H^0(\mathfrak g^1_6)\otimes \OO_S\rightarrow \mathfrak g^1_6$.
These bundles have played an important role in Voisin's recent proof
of Green's Conjecture on syzygies of canonical curves (cf.
\cite{V1}, \cite{V2}).
\end{remark}

The counterexample constructed in Theorem \ref{class} now raises at
least three questions:

\noindent $\bullet$ Is the divisor $\KK_{10}$ an isolated
counterexample? (After all, the condition that a curve lie on a $K3$
surface is divisorial only for $g=10$, and even on $\mm_{10}$ this
condition gives rise to a divisor almost by accident, due to the
somewhat miraculous existence of Mukai's rational $5$-fold
$X_5\subset \PP^{13}$).

\noindent $\bullet$ If the answer to the first question is no and
the Slope Conjecture fails systematically, are there integers $g\leq
23$ and divisors $D\in \mbox{Eff}(\mm_g)$ such that $s(D)<13/2$, so
that $\mm_g$ of general type, thus contradicting Conjecture
\ref{gentype}?

\noindent $\bullet$ In light of the application to the Schottky
problem, is there still a lower bound on $s(\mm_g)$? Note that we
know that $s(\mm_g)\geq \mbox{O}(1/g)$ for large $g$ (cf.
\cite{HMo}).

In the remaining part of this paper we will provide adequate answers
to the first two of these questions.

\section{Constructing divisors of small slope using syzygies}

We describe a general recipe  of constructing effective divisors on
$\mm_g$ having very small slope. In particular, we obtain an
infinite string of counterexamples to the Slope Conjecture.
Everything in this section is contained in \cite{F2} and \cite{F3}
and we outline the main ideas and steps in these calculations. The
key idea is to reinterpret the second description of the divisor
$\KK_{10}$ (see Theorem \ref{equiv}) as a failure of $[C]\in
\cM_{10}$ to satisfy the Green-Lazarsfeld property $(N_0)$ in the
embedding given by one of the finitely many linear series $\mathfrak
g^4_{12}$ on $C$. We will be looking at loci in $\cM_g$ consisting
of curves that have exceptional syzygy properties with respect to
certain $\mathfrak g^r_d$'s.

 Suppose that $C\stackrel{|L|}\hookrightarrow \PP^r$  is a curve of genus $g$ embedded by a
 line bundle $L\in \mbox{Pic}^d(C)$. We denote by $I_{C/\PP^r}$ the
ideal of $C$ in $\PP^r$ and consider its minimal resolution of free
graded $S=\mathbb C[x_0, \ldots, x_r]$-modules
$$0\rightarrow F_{r+1}\rightarrow \cdots \rightarrow F_2\rightarrow F_1\rightarrow
I_{C/\PP^r}\rightarrow 0.$$   Following Green and Lazarsfeld we say that the pair $(C, L)$
satisfy the property $(N_i)$ for some integer $i\geq 1$, if
$F_j=\oplus S(-j-1)$ for all $j\leq i$ (or equivalently in terms
of graded Betti numbers, $b_{i,l}(C)=0$ for all $l\geq 2$). Using
the computation of $b_{j,l}(C)$ in terms of Koszul cohomology,
there is a well-known cohomological interpretation of property
$(N_i)$: If $M_L$ is the vector bundle on $C$ defined by the exact sequence
$$ 0\rightarrow M_L\rightarrow H^0(L)\otimes \OO_C\rightarrow L\rightarrow
0,$$  then $(C,L)$ satisfies property $(N_i)$ if and only
if for all $j\geq 1$, the natural map
$$u_{i, j}:\wedge^{i+1} H^0(L)\otimes H^0(L^{\otimes j })\rightarrow H^0(\wedge ^i
M_L\otimes L^{\otimes (j+1)})$$
 is surjective (cf. e.g.
\cite{L2}).

Our intention is to define a determinantal syzygy type condition on
a generically finite cover of $\mm_g$ parametrizing pairs consisting
of a curve and a $\mathfrak g^r_d$. We fix integers $i\geq 0$ and
$s\geq 1$ and set $$r:=2s+si+i,\mbox{ } g:=rs+s \mbox{ and }
d:=rs+r.$$ We denote by $\mathfrak G^r_d$ the stack parametrizing
pairs $[C, L]$ with $[C]\in \cM_g$ and $L\in W^r_d(C)$ and denote by
$\sigma:\mathfrak G^r_d\rightarrow \cM_g$ the natural projection.
Since $\rho(g, r, d)=0$,  by general Brill-Noether theory, the
general curve of genus $g$ has finitely many $\mathfrak g^r_d$'s and
there exists a unique irreducible component of $\mathfrak G^r_d$
which maps onto $\cM_g$.

We define a  substack of $\mathfrak G^r_d$ consisting of those pairs
$(C, L)$ which fail to satisfy property $(N_i)$. In \cite{F3} we
introduced  two vector bundles $\cA$ and $\cB$ over $\mathfrak
G^r_d$ such that for a curve $C\stackrel{|L|}\hookrightarrow \PP^r$
corresponding to a point $(C, L)\in \mathfrak G^r_d$, we have that
$$\cA(C,L)=H^0(\PP^r, \wedge^i M_{\PP^r}(2)) \mbox{
}\mbox{ and }\cB(C,L)=H^0(C, \wedge^i M_L\otimes L^2).$$ There is a
natural vector bundle morphism $\phi:\cA\rightarrow \cB$ given by
restriction. From Grauert's Theorem we see that both $\cA$ and $\cB$
are vector bundles over $\mathfrak G^r_d$  and from Bott's Theorem
we compute their ranks
$$\mbox{rank}(\cA)=(i+1){r+2 \choose i+2}\mbox{  and  }\mbox{
rank}(\cB)={r\choose i}\Bigl(-\frac{id}{r}+2d+1-g\Bigr)$$ (use that
$M_L$ is a stable bundle, hence  $H^1(\wedge^i M_L\otimes L^{\otimes
2})=0$, while $\mbox{rank}(\cB)$ can be computed from Riemann-Roch).
It is easy to check that for our numerical choices we have that
$\mbox{rank}(\cA)=\mbox{rank}(\cB)$.

\begin{theorem}\label{ni}
The cycle $$\mathcal{U}_{g, i}:=\{(C,L)\in \mathfrak G^r_d
:(C,L) \mbox{ fails property }(N_i)\},$$ is the degeneracy locus of
vector bundle map $\phi:\cA\rightarrow \cB$ over $\mathfrak
G^r_d$.
\end{theorem}
Thus $\mathcal{Z}_{g, i}:=\sigma(\mathcal{U}_{g, i})$ is a virtual
divisor on $\cM_g$ when $g=s(2s+si+i+1)$. In \cite{F3} we show that
we can extend the determinantal structure of $\mathcal{Z}_{g, i}$
over $\mm_g$ in such a way that whenever $s\geq 2$, the resulting
virtual slope violates the Harris-Morrison Conjecture. One has the
following statement:

\begin{theorem}\label{div}
If $\sigma:\widetilde{\mathfrak G}^r_d\rightarrow \mm_g$ is the
compactification of $\mathfrak G^r_d$ given by limit linear series,
then there exists a natural extension of the vector bundle map
$\phi: \cA\rightarrow \cB$ over $\widetilde{\mathfrak G}^r_d$ such
that $\overline{\mathcal{Z}}_{g, i}$ is the image of the degeneracy
locus of $\phi$. The class of the pushforward to $\mm_g$ of the
virtual degeneracy locus of $\phi$ is given by
$$\sigma_*(c_1(\G_{i, 2}-\H_{i, 2}))\equiv a\lambda-b_0\ \delta_0-b_1\ \delta_1-\cdots -b_{[g/2]}\ \delta_{[g/2]},
  $$
where $a, b_0, \ldots, b_{[g/2]}$ are explicitly given coefficients
such that $b_1=12b_0-a$,  $b_i\geq b_0$ for $1\leq i\leq [g/2]$ and
$$s\bigl(\sigma_*(c_1(\G_{i, 2}-\H_{i, 2}))\bigr)=\frac{a}{b_0}=6\frac{f(s,
i)}{(i+2)\ s g(s, i)}, \mbox{ with} $$ \begin{center} $f(s,
i)=(i^4+24i^2+8i^3+32i+16)s^7+(i^4+4i^3-16i-16)s^6-(i^4+7i^3+13i^2-12)s^5-
(i^4+2i^3+i^2+14i+24)s^4
+(2i^3+2i^2-6i-4)s^3+(i^3+17i^2+50i+41)s^2+(7i^2+18i+9)s+2i+2$
\end{center}
and \noindent
\begin{center}
$g(s,i)=(i^3+6i^2+12i+8)s^6+(i^3+2i^2-4i-8)s^5-(i^3+7i^2+11i+2)s^4-\newline
-(i^3-5i)s^3+(4i^2+5i+1)s^2+ (i^2+7i+11)s+4i+2.$
\end{center}
Furthermore, we have that $6<\frac{a}{b_0}<6+\frac{12}{g+1}$
whenever $s\geq 2$. If the morphism $\phi$ is generically
non-degenerate, then $\overline{\mathcal{Z}}_{g, i}$ is a divisor on
$\mm_g$ which gives a counterexample to the Slope Conjecture for
$g=s(2s+si+i+1)$.
\end{theorem}

\begin{remark}
Theorem \ref{div} generalizes all known examples of effective
divisors on $\mm_g$ violating the Slope Conjecture. For $s=2$ and
$g=6i+10$ (that is, in the case $h^1(L)=2$ when $\mathfrak G^r_d$ is
isomorphic to a Hurwitz scheme parametrizing covers of $\PP^1$), we
recover our result from \cite{F2}. We have that
$$s(\overline{\mathcal{Z}}_{6i+10,
i})=\frac{3(4i+7)(6i^2+19i+12)}{(12i^2+31i+18)(i+2)}.$$
\end{remark}

For $i=0$ we recover the main result from \cite{Kh} originally proved using
a completely different method:
\begin{corollary}\label{khosla} (Khosla)
For $g=s(2s+1), r=2s, d=2s(s+1)$ the slope of the virtual class of
the locus of those $[C]\in \mm_g$ for which there exists $L\in
W^r_d(C)$ such that the embedded curve $C\subset \PP^r$ sits on a
quadric hypersurface, is
$$s(\overline{\mathcal{Z}}_{s(2s+1),
0})=\frac{3(16s^7-16s^6+12s^5-24s^4-4s^3+41s^2+9s+2)}{s
(8s^6-8s^5-2s^4+s^2+11s+2)}.$$
\end{corollary}

\begin{remark}
In the case $s=1, g=2i+3$ when $\mathfrak g^r_d=\mathfrak
g^{g-1}_{2g-2}$ is the canonical system, our formula reads
$$s(\overline{\mathcal{Z}}_{2i+3, i})=\frac{6(i+3)}{i+2}=6+\frac{12}{g+1}.$$
Remembering that $\mathcal{Z}_{2i+3, i}$ is the locus of curves
$[C]\in \cM_{2i+3}$ for which $K_C$ fails property $(N_i)$, from
Green's Conjecture for generic curves (cf. \cite{V1}, \cite{V2}) we
obtain the set-theoretic identification identification between
$\mathcal{Z}_{2i+3, i}$ and the locus $\mathcal{M}_{2i+3, i+2}^1$ of
$(i+2)$-gonal curves. Thus $\mathcal{Z}_{2i+3, i}$ is a
Brill-Noether divisor! Theorem \ref{div} provides a new way of
calculating the class of the compactification of the Brill-Noether
divisor which was first computed by Harris and Mumford (cf.
\cite{HM}).
\end{remark}

Theorem \ref{div} is proved by extending  the determinantal
structure of $\mathcal{Z}_{g, i}$ over the boundary divisors in
$\mm_{g}$. We can carry this out outside a locus of codimension
$\geq 2$ in $\mm_g$. We denote by $\widetilde{\cM}_g:=\cM_g^0\cup
\bigl(\cup_{j=0}^{[g/2]} \Delta_j^0\bigr)$ the locally closed subset
of $\mm_g$ defined as the union  of the locus $\cM_{g}^0$ of smooth
curves carrying no linear systems $\mathfrak g^r_{d-1}$ or
$\mathfrak g^{r+1}_d$ to which we add the open subsets
$\Delta_j^0\subset \Delta_j$ for $1\leq j\leq [g/2]$ consisting of
$1$-nodal genus $g$ curves $C\cup_y D$, with $[C]\in \cM_{g-j}$ and
$[D, y]\in \cM_{j, 1}$ being Brill-Noether general curves, and the
locus $\Delta_0^0\subset \Delta_0$ containing $1$-nodal irreducible
genus $g$ curves $C'=C/q\sim y$, where $[C, q]\in \cM_{g-1}$ is a
Brill-Noether general pointed curve and $y\in C$, together with
their degenerations consisting of unions of a smooth genus $g-1$
curve and a nodal rational curve. One can then extend the finite
covering $\sigma:\mathfrak G^{r}_{d}\rightarrow \cM_g^0$ to a
proper, generically finite map
$$\sigma: \widetilde{\mathfrak G}^{r}_{d} \rightarrow \widetilde{\cM}_g$$ by
letting $\widetilde{\mathfrak G}^{r}_{d}$ be the stack of
limit $\mathfrak g^{r}_{d}$'s on the treelike curves from
$\widetilde{\cM}_g$ (see \cite{EH1}, Theorem 3.4 for the
construction of the space of limit linear series).

One method of computing $[\overline{\mathcal{Z}}_{g, i}]$ is to
intersect the locus $\overline{\mathcal{Z}}_{g, i}$ with standard
test curves in the boundary of $\mm_{g}$ which are defined as
follows: we fix a Brill-Noether general curve $C$ of genus $g-1$, a
general point $q\in C$ and a general elliptic curve $E$. We define
two $1$-parameter families
\begin{equation}\label{testcurves}
C^0:=\{C/y\sim q: y\in C\}\subset \Delta_0 \subset \mm_{g}
\mbox{ and }C^1:=\{C\cup _y E: y\in C\}\subset \Delta_1\subset
\mm_{g}.
\end{equation}
 It is well-known that
these families intersect the generators of
$\mbox{Pic}(\mm_{g})$ as follows:
$$
C^0\cdot \lambda=0,\ C^0\cdot \delta_0=-(2g-2), \ C^0\cdot
\delta_1=1 \mbox{ and } C^0\cdot \delta_a=0\mbox{ for }a\geq 2,
\mbox{ and}$$
$$C^1\cdot \lambda=0, \ C^1\cdot \delta_0=0, \ C^1\cdot
\delta_1=-(2g-4), \ C^1\cdot \delta_a=0 \mbox{ for }a\geq 2.$$

 Before we proceed we review the notation used in the theory of limit linear series
 (see \cite{EH1} as a general reference).
 If $X$ is a
treelike curve and $l$ is a limit $\mathfrak g^r_d$ on $X$, for a
component  $Y$ of $X$ we denote by $l_Y=(L_Y, V_Y\subset H^0(L_Y))$
the $Y$-aspect of $l$. For a point $y\in Y$ we denote by by
$\{a^{l_Y}_s(C)\}_{s=0\ldots r}$ the \emph{vanishing sequence} of
$l$ at $y$ and by $\rho(l_Y, y):=\rho(g, r, d)-\sum_{i=0}^r
(a^{l_Y}_i(y)-i)$ the adjusted Brill-Noether number with respect to
$y$. We have the following description of the curves $\sigma^*(C^0)$
and $\sigma^*(C^1)$:
\begin{proposition}\label{limitlin}
(1) Let $C_y^1=C\cup_y E$ be an element of $\Delta_1^0$. If $(l_C,
l_E)$ is a limit  $\mathfrak g^{r}_{d}$ on $C_y^1$, then
$V_C=H^0(L_C)$ and $L_C\in W^{r}_{d}(C)$ has a cusp at $y$. If $y\in
C$ is a general point, then $l_E=\bigl(\OO_E(dy),
(d-r-1)y+|(r+1)y|\bigr)$, that is, $l_E$ is uniquely determined. If
$y\in C$ is one of the finitely many points for which there exists
$L_C\in W^{r}_{d}(C)$ such that $\rho(L_C, y)=-1$, then
$l_E(-(d-r-2)y)$ is a $\mathfrak g_{r+2}^{r}$ with vanishing
sequence at $y$ being $\geq (0, 2, 3, \ldots, r, r+2)$.  Moreover,
at the level of $1$-cycles we have the identification
$\sigma^*(C^1)\equiv X + \nu \ T$, where
$$X:=\{(y, L)\in C\times W^r_d(C) :h^0(C, L(-2y))\geq r\}$$
and $T$ is the curve consisting of $\mathfrak g^{r}_{r+2}$'s on $E$
with vanishing $\geq (0, 2, \ldots, r, r+2)$ at the fixed point
$y\in E$ while $\nu$ is a positive integer.

\noindent (2) Let $C_y^0=C/y\sim q$ be an element of $\Delta_0^0$.
Then limit linear series of type $\mathfrak g^{r}_{d}$ on $C_y^0$
are in 1:1 correspondence with complete linear series $L$ on $C$ of
type $\mathfrak g^{r}_{d}$ satisfying the condition $h^0(C, L\otimes
\OO_C(-y-q))=h^0(C,L)-1.$ Thus there is an isomorphism between the
cycle $\sigma^*(C^0)$ of $\mathfrak g^{r}_{d}$'s on all curves
$C_y^0$ with $y\in C$, and the smooth curve
$$Y:=\{(y, L)\in C\times W^r_d(C): h^0(C, L(-y-q))\geq r\}.$$

\end{proposition}

Throughout the papers \cite{F2} and \cite{F3} we use a number of
facts about intersection theory on Jacobians which we now quickly
review. Let $C$ be a Brill-Noether general curve of genus $g-1$
(recall that $g=rs+s$ and $d=rs+s$, where $r=2s+si+i$). Then
$\mbox{dim } W^r_d(C)=r$ and every $L\in W^r_d(C)$ corresponds to a
complete and base point free linear series.  We denote by $\L$ a
Poincar\'e bundle on $C\times \mbox{Pic}^d(C)$ and by $\pi_1:C\times
\mbox{Pic}^d(C)\rightarrow C$ and $\pi_2:C\times
\mbox{Pic}^d(C)\rightarrow \mbox{Pic}^d(C)$ the projections. We
define the cohomology class $\eta=\pi_1^*([\mbox{point}])\in
H^2(C\times \mbox{Pic}^d(C))$, and if $\delta_1,\ldots,
\delta_{2g}\in H^1(C, \mathbb Z)\cong H^1(\mbox{Pic}^d(C), \mathbb
Z)$ is a symplectic basis, then we set
$$\gamma:=-\sum_{\alpha=1}^g
\Bigl(\pi_1^*(\delta_{\alpha})\pi_2^*(\delta_{g+\alpha})-\pi_1^*(\delta_{g+\alpha})\pi_2^*(\delta_
{\alpha})\Bigr).$$ We have the formula (cf. \cite{ACGH}, p. 335)
$c_1(\L)=d\eta+\gamma,$ corresponding to the Hodge decomposition of
$c_1(\L)$. We also record that $\gamma^3=\gamma \eta=0$, $\eta^2=0$
and $\gamma^2=-2\eta \pi_2^*(\theta)$. On $W^r_d(C)$  we have the
tautological rank $r+1$ vector bundle
$\mathcal{E}:=(\pi_2)_{*}(\mathcal{L}_{| C\times W^r_d(C)})$. The
Chern numbers of $\mathcal{E}$ can be computed using the Harris-Tu
formula (cf. \cite{HT}): if we write $\sum_{i=0}^r
c_i(\mathcal{E}^{\vee})=(1+x_1)\cdots (1+x_{r+1})$, then for every
class $\zeta \in H^*(\mbox{Pic}^d(C), \mathbb Z)$ one has the formula
\footnote{There is a confusing sign error in the formula (1.4) in
\cite{HT}: the formula is correct as it is appears in \cite{HT}, if
the $x_j$'s denote the Chern roots of the \emph{dual} of the kernel
bundle.}
$$x_1^{i_1}\cdots x_{r+1}^{i_{r+1}}\
\zeta=\mbox{det}\Bigl(\frac{\theta^{g-1+r-d+i_j-j+l}}{(g-1+r-d+i_j-j+l)!}\Bigr)_{1\leq
j, l\leq r+1}\ \zeta.$$ If  we use the expression of the Vandermonde
determinant, we get the formula
$$\mbox{det}\Bigl(\frac{1}{(a_j+l-1)!}\Bigr)_{1\leq j, l\leq
r+1}=\frac{\Pi_{ j>l}\ (a_l-a_j)}{\Pi_{j=1}^{r+1}\ (a_j+r)!}.$$ By
repeatedly applying this we get all intersection numbers on
$W^r_d(C)$ which we need:
\begin{lemma}\label{vandermonde}
If $c_i:=c_i(\mathcal{E}^{\vee})$ we have the following identities
in $H^*(W^r_d(C), \mathbb Z)$:
\begin{enumerate}
\item $c_{r-1}\theta=\frac{r(s+1)}{2} c_r.$
\item $c_{r-2} \theta^2=\frac{r(r-1)(s+1)(s+2)}{6} c_r.$
\item $c_{r-2} c_1 \theta=\frac{r(s+1)}{2}\bigl(1+\frac{(r-2)(r+2)(s+2)}{3(s+r+1)}\bigr)
c_r$.
\item $c_{r-1} c_1=(1+\frac{(r-1)(r+2)(s+1)}{2(s+r+1)})c_r.$
\item $c_r=\frac{1! \ 2!\cdots (r-1)!\ (r+1)!}{(s-1)!\  (s+1)!\ (s+2)! \cdots
(s+r)!}\theta^{g-1}.$
\end{enumerate}
\end{lemma}

For each integers $0\leq a\leq r$ and $ b\geq 2$ we shall define
vector bundles $\G_{a,b}$ and $\H_{a, b}$ over $\widetilde{\mathfrak
G}^{r}_{d}$ with fibres
$$\G_{a,b}(C,L)=H^0(C, \wedge^a M_L\otimes L^{\otimes b})\ \mbox{ and }\ \H_{a, b}(C,L)=H^0(\PP^r,
\wedge^a M_{\PP^r}(b))$$ for each $(C, L)\in \mathfrak G^{r}_{d}$
giving a map $C\stackrel{|L|} \rightarrow \PP^r$.  Clearly $\G_{i, 2
| \mathfrak G^{r}_{d}}=\mathcal{B}$ and $\H_{i, 2 | \mathfrak
G^r_d}=\cA$, where $\cA$ and $\mathcal{B}$ are the vector bundles
introduced in Proposition \ref{ni}. The question is how to extend
this description over the divisors $\Delta_j^0$. For simplicity we
only explain how to do this over $\sigma^{-1}(\cM_g^0\cup
\Delta_0^0\cup \Delta_1^0)$ which will be enough to compute the
slope of $\overline{\mathcal{Z}}_{g,i}$. The case of the divisors
$\sigma^{-1}(\Delta_j^0)$ where $2\leq j\leq [g/2]$ is technically
more involved and it is dealt with in \cite{F3}.
 We start by extending $\G_{0,
b}$ (see \cite{F3}, Proposition 2.8):

\begin{proposition} For each $b\geq 2$ there
exists a vector bundle $\G_{0,b}$ over $\widetilde{\mathfrak
G}^{r}_{d}$ of rank $bd+1-g$ whose fibres admit the following
description:
\begin{itemize}
\item For $(C, L)\in \mathfrak G^{r}_{d}$, we have that
$\G_{0,b}\bigl(C,L)=H^0(C, L^{\otimes b})$.
\item For $t=(C\cup_y
E, L)\in \sigma^{-1}(\Delta_1^0)$, where $L\in W^r_d(C)$ has a cusp
at $y\in C$, we have that
$$\G_{0,b}(t)=H^0(C, L^{\otimes b}(-2y))+\mathbb C\cdot
u^b\subset H^0\bigr(C, L^{\otimes b}),$$ where $u\in H^0(C, L)$ is
any section such that $\rm{ord}$$_y(u)=0$. \item For $t=(C/y\sim q,
L)\in \sigma^{-1}(\Delta_0^0)$, where $q,y\in C$ and $L\in W^r_d(C)$
is such that $h^0(C, L(-y-q))=h^0(L)-1$, we have that
$$\G_{0,b}(t)=H^0(C, L^{\otimes b}(-y-q))\oplus \mathbb C\cdot u^b\subset H^0(C, L^{\otimes b}),$$
where $u\in H^0(C, L)$ is a section such that $\rm{ord}$$_y(u)=\rm{ord}$$_q(u)=0$.
\end{itemize}
\end{proposition}

Having defined the vector bundles $\G_{0, b}$ we now define
inductively all vector bundles $\G_{a, b}$ by the exact sequence
\begin{equation}\label{gi}
0\longrightarrow \G_{a, b}\longrightarrow \wedge^a \G_{0,
1}\otimes \G_{0, b}\stackrel{d_{a, b}}\longrightarrow \G_{a-1,
b+1}\longrightarrow 0.
\end{equation}
To define $\H_{a, b}$ is even easier. We set $\H_{0,
b}:=\mbox{Sym}^b \G_{0, 1}$ for all $b\geq 1$ and we define $\H_{a,
b}$ inductively via the exact sequence
\begin{equation}\label{hi}
0\longrightarrow \H_{a, b}\longrightarrow \wedge^a \H_{0, 1}\otimes
\mbox{Sym}^b \H_{0, 1}\longrightarrow \H_{a-1, b+1}\longrightarrow
0.
\end{equation}

The surjectivity of the right map in (\ref{hi}) is obvious, whereas
to prove that $d_{a, b}$ is surjective, one employs the arguments
from Proposition 3.10 in \cite{F2}. There is a natural vector bundle
morphism $\phi_{a, b}:\H_{a, b}\rightarrow \G_{a, b}$. Moreover
$\mbox{rank}(\H_{i, 2})=\mbox{rank}(\G_{i, 2})$ and the degeneracy
locus of $\phi_{i, 2}$ is the codimension one compactification of
$\mathcal{Z}_{g, i}$.

We now compute the class of the curves $X$ and $Y$ defined in
Proposition \ref{limitlin} (see \cite{F3} Proposition 2.11 for
details):

\begin{proposition}\label{xy}
Let $C$ be a Brill-Noether general curve of genus $g-1$ and $q\in C$
a general point. We denote by $\pi_2: C\times W^r_d(C)\rightarrow
W^r_d(C)$ the projection and set
$c_i:=(\pi_2)^*\bigl(c_i(\mathcal{E}^{\vee})\bigr)$.

\noindent(1) The class of the curve $X=\{(y,L)\in C\times
W^r_d(C):h^0(C,L(-2y))\geq r\}$ is given by
$$[X]= c_r+c_{r-1}(2\gamma+(2d+2g-4)\eta)-6c_{r-2}\  \eta \theta.
$$

\noindent (2) The class of the curve $Y=\{(y,L)\in C\times
W^r_d(C):h^0(C,L(-y-q))\geq r\}$ is given by
$$[Y]= c_r+c_{r-1}(\gamma +(d-1)\eta)-2c_{r-2}\ \eta \theta.$$
\end{proposition}
\noindent{\emph{Sketch of proof.}} Both $X$ and $Y$ are expressed as
degeneracy loci over $C\times W^r_d(C)$ and we compute their classes
using the Thom-Porteous formula. For $(y, L)\in C\times W^r_d(C)$
the natural map $ H^0(C, L_{|2y})^{\vee} \rightarrow H^0(C,
L)^{\vee}$ globalizes to a vector bundle map $\zeta:
J_1(\mathcal{L})^{\vee} \rightarrow (\pi_2)^*(\mathcal{E}^{\vee})$.
Then $X=Z_1(\zeta)$ and we apply Thom-Porteous. \hfill $\Box$

We mention the following intersection theoretic result (cf.
\cite{F3}, Lemma 2.12):
\begin{lemma}\label{fj}
For each $j\geq 2$ we have the following \texttt{}formulas:
\begin{enumerate}
\item
$c_1(\G_{0, j \ |X})=-j^2\theta-(2g-4)\eta-j(d\eta+\gamma)$.
\item
$c_1(\G_{0, j \ |Y})=-j^2\theta+\eta$.
\end{enumerate}
\end{lemma}

\noindent \emph{Proof of Theorem \ref{div}.}
 Since
$\mbox{codim}(\mm_g-\widetilde{\cM}_g, \mm_g)\geq 2$, it makes no
difference whether we compute the class $\sigma_*(\G_{i, 2}-\H_{i,
2})$ on $\widetilde{\cM}_g$ or on $\mm_g$ and we can write
$$\sigma_*(\G_{i, 2}-\H_{i, 2})= A \lambda-B_0\ \delta_0-B_1\ \delta_1-\cdots
-B_{[g/2]}\ \delta_{[g/2]}\in \mbox{Pic}(\mm_g),$$
 where $\lambda, \delta_0, \ldots, \delta_{[g/2]}$ are the
generators of $\mbox{Pic}(\mm_g)$. First we note that one has the
relation $A-12B_0+B_1=0$. This can be seen by picking a general
curve $[C, q]\in \cM_{g-1, 1}$ and at the fixed point $q$  attaching
to $C$ a Lefschetz pencil of plane cubics. If we denote by $R\subset
\mm_g$ the resulting curve, then
 $R\cdot \lambda=1, \ R\cdot \delta_0=12, \ R\cdot
\delta_1=-1$ and $R\cdot \delta_j=0$ for $j\geq 2$. The relation
$A-12B_0+B_1=0$ follows once we show that $\sigma^*(R)\cdot
c_1(\G_{i, 2}-\H_{i, 2})=0$. To achieve this we check that $\G_{0, b
|\sigma^*(R)}$ is trivial and then use (\ref{gi}) and (\ref{hi}). We
take $[C\cup_q E]$ to be an arbitrary curve from $R$, where $E$ is
an elliptic curve. Using that limit $\mathfrak g^{r}_{d}$ on
$C\cup_q E$ are in $1:1$ correspondence with linear series $L\in
W^{r}_{d}(C)$ having a cusp at $q$ (this being a statement that
holds independent of $E$) and that $\G_{0, b |
\sigma^*(\Delta_1^0)}$ consists on each fibre of sections of the
genus $g-1$ aspect of the limit $\mathfrak g^{r}_{d}$, the claim now
follows.

Next we  determine $A, B_0$ and $B_1$ explicitly.  We fix a general
pointed curve $(C, q)\in \cM_{g-1, 1}$ and construct the test curves
$C^1\subset \Delta_1$ and $C^0\subset \Delta_0$. Using the notation
from Proposition \ref{limitlin}, we get that $\sigma^*(C^0)\cdot
c_1(\G_{i, 2}-\H_{i, 2}) =c_1(\G_{i, 2 |Y})-c_1(\H_{i, 2 |Y})$   and
 $\sigma^*(C^1) \cdot c_1(\G_{i, 2}-\H_{i, 2})=c_1(\G_{i, 2 |X})-c_1(\H_{i, 2 |X})$
 (the other component $T$ of $\sigma^*(C^1)$ does not appear
 because $\G_{0, b |T}$ is trivial for all $b\geq 1$).
On the other hand $$C^0 \cdot \sigma_*(c_1(\G_{i, 2}-\H_{i,
2}))=(2g-2)B_0-B_1 \mbox{ and } C^1\cdot \sigma_*(c_1(\G_{i,
2}-\H_{i, 2}))=(2g-4)B_1,$$ while we already know that
$A-12B_0+B_1=0.$ Next we use the relations
$$c_1(\G_{i, 2})= \sum_{l=0}^i (-1)^l c_1(\wedge^{i-l}\G_{0,
1}\otimes \G_{0, l+2})=\sum_{l=0}^i (-1)^l{ r+1 \choose i-l}
c_1(\G_{0, l+2})+$$ $$+\sum_{l=0}^i (-1)^l
\bigl((l+2)(rs+r)+1-rs-s\bigr){r \choose i-l-1}c_1(\G_{0, 1}), \
\mbox{ } \mbox{ and }$$
$$c_1(\H_{i, 2})=\sum_{l=0}^i (-1)^l c_1(\wedge^{i-l} \G_{0,
1}\otimes \mbox{Sym}^{l+2} \G_{0, 1})=$$ $$=\sum_{l=0}^i
(-1)^l\Bigl( {r \choose i-l-1}{r+l+2 \choose l+2}+{r+1\choose
i-l}{r+l+2\choose r+1}\Bigr)c_1(\G_{0, 1}),$$ which when restricted
to $X$ and $Y$, enable us (also using Lemma \ref{fj}) to obtain
explicit expressions for $c_1(\G_{i, 2}-\H_{i, 2})_{|X}$ and
$c_1(\G_{i, 2}-\H_{i, 2})_{| Y}$ in terms of the classes $\eta,
\theta, \gamma$ and $c_1=\pi_2^*(c_1(\mathcal{E}^{\vee}))$.
Intersecting these classes with $[X]$ and $[Y]$ and using Lemma
\ref{vandermonde}, we finally get a linear system of $3$ equations
in $A, B_0$ and $B_1$ which leads to the stated formulas for the
first three coefficients. \hfill $\Box$

Theorem \ref{class} produces only virtual divisors on $\mm_g$ of
slope $<6+12/(g+1)$. To get actual divisors one has to show  that
the vector bundle map $\phi:\H_{i, 2}\rightarrow \G_{i, 2}$ is
generically non-degenerate. This has been carried out for $s=2, i=0$
in \cite{FP} (relying on earlier work by Mukai), as well in the
cases $s=2, i=1$ and $s=2, i=2$ in \cite{F2}, using the program
Macaulay. D. Khosla has checked the transversality of $\phi$ when
$s=3, i=0$, that is on $\cM_{21}$ (cf. \cite{Kh}). We generalize
this last result as well as \cite{FP} by proving that for $i=0$ and
arbitrary $s$, the map $\phi: \H_{0, 2}\rightarrow \G_{0, 2}$ is
always generically non-degenerate:

\begin{theorem}\label{maxrank}
For an integer $s\geq 2$ we set $r:=2s, d:=2s(s+1)$ and
$g:=s(2s+1)$. Then the vector bundle map $\phi:\H_{0, 2}\rightarrow
\G_{0, 2}$ is generically non-degenerate. In particular
$$\mathcal{Z}_{g, 0}:=\{[C]\in \cM_g: \exists L\in W^r_d(C) \mbox{ such that } C\stackrel{|L|}\hookrightarrow \PP^r
\mbox{ is not projectively normal}\}$$ is a divisor on $\mm_g$ of
slope
$$s(\overline{\mathcal{Z}}_{g,
0})=\frac{3(16s^7-16s^6+12s^5-24s^4-4s^3+41s^2+9s+2)}{s(8s^6-8s^5-2s^4+s^2+11s+2)}$$
violating the Slope Conjecture.
\end{theorem}
\noindent \emph{Sketch of proof.} From Brill-Noether theory it
follows that there exists a unique component of
$\widetilde{\mathfrak G}^r_d$ which maps onto $\widetilde{\cM}_g$,
therefore it is then enough to produce a Brill-Noether-Petri general
smooth curve $C\subset \PP^{2s}$ having degree $2s(s+1)$ and genus
$s(2s+1)$ such that $C$ does not sit on any quadrics, that is
$H^0(\mathcal{I}_{C/\PP^{2s}}(2))=H^1(\mathcal{I}_{C/\PP^{2s}}(2))=0$.
We carry this out inductively: for each $0\leq a\leq s$, we
construct a smooth non-degenerate curve $C_a \subset \PP^{s+a}$ with
$\mbox{deg}(C_a)={s+a+1 \choose 2}+a$ and $g(C_a)={s+a+1\choose
2}+a-s$, such that $C_a$ satisfies the Petri Theorem (in particular
$H^1(C_a, N_{C_a /\PP^{s+a}})=0$), and such that the multiplication
map $\mu_2: \mbox{Sym}^2 H^0(C_a, \OO_{C_a}(1))\rightarrow H^0(C_a,
\OO_{C_a}(2))$ is surjective.

To construct $C_0\subset \PP^s$ we consider the \emph{White surface}
$S=\mbox{Bl}_{\{p_1, \ldots, p_{\delta}\}}(\PP^2) \subset \PP^s$
obtained by blowing-up $\PP^2$ at general points $p_1, \ldots,
p_{\delta}\in \PP^2$ where $\delta={s+1\choose 2}$,  and embedding
it via the linear system $|s h-\sum_{i=1}^{\delta} E_{p_i}|$. Here
$h$ is the class of a line on $\PP^2$. It is known that $S\subset
\PP^s$ is a projectively Cohen-Macaulay surface and its ideal is
generated by the $(3\times 3)$-minors of a certain $(3\times
s)$-matrix of linear forms. The Betti diagram of $S\subset \PP^s$ is
the same as that of the ideal of $(3\times 3)$-minors of a $(3\times
s)$-matrix of indeterminates. In particular, we have that
$H^i(\mathcal{I}_{S/\PP^s}(2))=0$ for $i=0, 1$. On $S$ we consider a
generic smooth curve $C \equiv (s+1)h-\sum_{i=1}^{\delta} E_{p_i}$.
We find that the embedded curve $C\subset S\subset \PP^s$ has
$\mbox{deg}(C)={s+1\choose 2}$ and $g(C)={s \choose 2}$.
 Even though $[C]\in \cM_{g(C)}$ itself is not a Petri
general curve, the map $H_{d(C), g(C), s}\rightarrow \cM_{{s\choose
2}}$ from the Hilbert scheme of curves
 $C'\subset \PP^s$, is smooth and dominant around the point $[C\hookrightarrow
 \PP^s]$. Therefore a generic deformation $[C_0\hookrightarrow
 \PP^s]$ of $[C\hookrightarrow \PP^s]$ will be Petri general and still satisfy
the condition $H^1(\mathcal{I}_{C_0/\PP^s}(2))=0$.

Assume now that we have constructed a Petri general curve
$C_a\subset \PP^{s+a}$ with all the desired properties. We pick
general points $p_1, \ldots, p_{s+a+2}\in C_a$ with the property
that if $\Delta:=p_1+\cdots +p_{s+a+2}\in \mbox{Sym}^{s+a+2}C_{a}$,
then the variety
$$T:=\{(M, V)\in W^{s+a+1}_{d(C_a)+s+a+2}(C_a):
\mbox{dim}\bigl(V\cap H^0(C_a, M\otimes
\OO_{C_a}(-\Delta))\bigr)\geq s+a+1\}$$ of linear series having an
$(s+a+2)$-fold point along $\Delta$, has the expected dimension
$\rho(g(C_a), s+a+1, d(C_a)+s+a+2)-(s+a+1)^2$. We identify the
projective space $\PP^{s+a}$ containing $C_a$ with a hyperplane
$H\subset \PP^{s+a+1}$ and choose a linearly normal elliptic curve
$E\subset \PP^{s+a+1}$ such that $E\cap H=\{p_1, \ldots,
p_{s+a+2}\}$. We set $X:=C_a\cup_{\{p_1, \ldots,
p_{s+a+2}\}}E\hookrightarrow \PP^{s+a+1}$ and then
$\mbox{deg}(X)=p_a(X)+s$. From the exact sequence $$0\longrightarrow
\OO_E(-p_1-\cdots -p_{s+a+2})\longrightarrow \OO_X\longrightarrow
\OO_{C_a}\longrightarrow 0,$$ we can write that $h^0(\OO_X(1))\leq
h^0(\OO_{C_a}(1))+h^0(\OO_E)=s+a+2$, hence $h^0(\OO_X(1))=s+a+2$ and
$h^1(\OO_X(1))=a+1$. One can also write the exact sequence
$$0\longrightarrow \mathcal{I}_{E/\PP^{s+a+1}}(1)\longrightarrow
\mathcal{I}_{X/\PP^{s+a+1}}(2)\longrightarrow
\mathcal{I}_{C_a/H}(2)\longrightarrow 0,$$ from which we obtain that
$H^1(\mathcal{I}_{X/\PP^{s+a+1}}(2))=0$, hence by a dimension count
also $H^0(\mathcal{I}_{X/\PP^{s+a+1}}(2))=0$, that is, $X$ and every
deformation of $X$ inside $\PP^{s+a+1}$ will lie on no quadrics. In
\cite{F3} Theorem 1.5 it is proved that $X\hookrightarrow
\PP^{s+a+1}$ can be deformed to an embedding of a smooth curve
$C_{a+1}$ in $\PP^{s+a+1}$ such that
$H^1(N_{C_{a+1}/\PP^{s+a+1}})=0$. This enables us to continue the
induction and finish the proof.\hfill $\Box$

\section{The Kodaira dimension of $\mm_g$ and other problems}

Since one is able to produce systematically effective divisors on
$\mm_g$ having slope
 smaller than that of the Brill-Noether divisors, it is natural to ask whether one could diprove Conjecture
 \ref{gentype}, that is, construct
effective divisors $D\in \mbox{Eff}(\mm_g)$ for $g\leq 23$ such that $s(D)<s(K_{\mm_g})=13/2$,
which would imply that $\mm_g$ is of general type. We almost succeeded in this with Theorem \ref{div} in
 the case $s=2, i=2, g=22$: The slope of the (actual) divisor $\overline{\mathcal{Z}}_{22, 2}\subset \mm_{22}$ turns out to
 be $1665/256=6.5039...,$ which barely fails to make $\mm_{22}$ of general type. However, a different syzygy
 type condition, this time pushed-forward from a variety which maps onto $\mm_{22}$ with fibres of dimesnion
 one, produces an effective divisor of slope even smaller
than $s(\overline{\mathcal{Z}}_{22, 2})$. We have the following
result \cite{F4}:

\begin{theorem}
The moduli space $\mm_{22}$ is of general type. Precisely, the locus
$$D_{22}:=\{[C]\in \cM_{22}: \exists L\in W^6_{25}(C) \mbox{ such that } C\stackrel{|L|}\hookrightarrow \PP^6
\mbox{ lies on a quadric}\}$$ is a divisor on $\cM_{22}$ and the
class of its closure in $\mm_{22}$ equals
$$\overline{D}\equiv c(\frac{17121}{2636}\lambda-\delta_0-\frac{14511}{2636}\delta_1-b_2\ \delta_2-\cdots
-b_{11}\ \delta_{11}),$$
where $c>0$ and $b_i>1$ for $2 \leq i \leq 11$. Therefore $s(\overline{D})=17121/2636=6.49506...<13/2.$
\end{theorem}

We certainly expect a similar result for $\mm_{23}$. We have
calculated the class of the virtual locus $D_{23}$ consisting of
curves $[C]\in \cM_{23}$ such that there exists $L\in W^6_{26}(C)$
with the multiplication map $\mu_L: \mbox{Sym}^2 H^0(L)\rightarrow
H^0(L^{\otimes 2})$ not being injective.  By dimension count we
expect this locus to be a divisor on $\mm_{23}$ and assuming so, we
have computed its slope $s(\overline{D}_{23})=
470749/72725=6.47300...<13/2$. For $g=23$ this is only a virtual
result at the moment, since we cannot rule out the possibility that
$D_{23}$ equals the entire moduli space $\cM_{23}$. The difficulty
lies in the fact that $D_{23}$ as a determinantal variety is
expected to be of codimension $3$ inside the  variety
$\mathfrak{G}^6_{26}$ which maps onto $\cM_{23}$ with fibres of
dimension $2$.

\noindent{\bf{The Kodaira dimension of $\mm_{g, n}$}} \vskip 4pt
\noindent The problem of describing the Kodaira type of $\mm_{g, n}$
for $n\geq 1$, has been initiated by Logan in \cite{Log}. Using
Theorem \ref{canonical}  together with the formula $K_{\mm_{g,
n}}=\pi_n^*(K_{\mm_{g, n-1}})+\omega_{\pi_n}$, where $\pi_n:\mm_{g,
n}\rightarrow \mm_{g, n-1}$ is the projection map forgetting the
$n$-th marked point, we find that
$$K_{\mm_{g, n}}\equiv 13\lambda-2\delta_0+\sum_{i=1}^n \psi_i-2\sum_{i\geq 0, S} \delta_{i: S}-\sum_{S}\delta_{1:S}.$$
To prove that $\mm_{g, n}$ is of general type one needs an ample
supply of explicit effective divisor classes on $\mm_{g, n}$ such
that $K_{\mm_{g, n}}$ can be expressed as a linear combination with
positive coefficients of such an effective divisor, boundary classes
and an ample class on $\mm_{g, n}$ of the type $\sum_{i=1}^n a_i
\psi_i+b \lambda-\delta_0-\sum_{i\geq 0, S} \delta_{i: S}$ where
$b>11$ and $a_i>0$ for $1\leq i\leq n$. Logan has computed the class
of the following effective divisors on $\mm_{g, n}$ (cf. \cite{Log},
Theorems 5.3-5.7): We fix nonnegative integers $a_1, \ldots, a_n$
such that $a_1+\cdots+a_n=g$ and we define $D_{g: a_1, \ldots, a_n}$
to be the locus of curves $[C, x_1, \ldots, x_n]\in \cM_{g, n}$ such
that $h^0(C, \OO_C(a_1 x_1+\cdots +a_n x_n))\geq 2$. Then $D_{g:
a_1, \ldots, a_n}$ is a divisor on $\mm_{g, n}$ and one has the
following formula in $\mbox{Pic}(\mm_{g, n})$:
\begin{equation}\label{logan} \overline{D}_{g: a_1, \ldots,
a_n}\equiv -\lambda+\sum_{i+1}^n {a_i+1\choose 2}\psi_i-0\cdot
\delta_0- \sum_{i<j} {a_i+a_j+1\choose 2} \delta_{0:\{i,
j\}}-\cdots.
\end{equation}
 Note that for $n=1$ when necessarily $a_1=g$, we obtain in
this way the class of the divisor of Weierstrass points on $\mm_{g,
1}$: $\overline{D}_{g:g}\equiv -\lambda+{g+1\choose 2}-\sum_{i=1}^g
{g-i+1\choose 2}\delta_{i:1}. $

In \cite{F3} we introduced a new class of divisors generalizing the
loci of higher Weierstrass points in a different way: Fix $g, r\geq
1$ and $0\leq i\leq g$. We set $n:=(2r+1)(g-1)-2i$ and define the
locus
$$\mathfrak{Mrc}_{g,i}^r:=\{[C,x_1, \ldots, x_n]\in \cM_{g,n}: h^1\bigl(C, \wedge^{i}
M_{K_C}\otimes K_C^{\otimes (r+1)}\otimes
\OO_C(-x_1-\cdots-x_n)\bigr)\geq 1\}.$$

If we denote by $\Gamma:=x_1+\cdots+x_n\in C_n$, by Serre duality,
the condition appearing in the definition of $\mathfrak{Mrc}_{g,
i}^r$ is equivalent to
$$
h^0\bigl(C, \wedge^{i} M_{K_C}^{\vee}\otimes
\OO_C(\Gamma)\otimes K_C^{\otimes (-r)}\bigr)\geq 1
\Longleftrightarrow \OO_C(\Gamma)\otimes K_C^{\otimes (-r)}\in
\Theta_{\wedge^{i} M_{K_C}^{\vee}},$$ where we recall that for a
stable vector bundle $E$ on $C$ having  slope $\nu(E)=\nu\in \mathbb
Z$, its \emph{theta divisor} is the determinantal locus
$$\Theta_E:=\{ \eta \in \mbox{Pic}^{g-\mu-1}(C): h^0(C, E\otimes
\eta) \geq 1 \}.$$ The main result from \cite{FMP} gives an
identification $\Theta_{\wedge^i M_{K_C}^{\vee}}=C_{g-i-1}-C_i$,
where the right hand side is one of the difference varieties
associated to $C$. Thus one has an alternative description of points
in $\mathfrak{Mrc}_{g,i}^r$: a point $(C, x_1, \ldots, x_n)\in
\mathfrak{Mrc}_{g, i}^r$ if and only if there exists $D\in C_i$ such
that $h^0\bigl(C, \OO_C(\Gamma+D)\otimes K_C^{\otimes
(-r)}\bigr)\geq 1$. For $i=0$, the divisor $\mathfrak{Mrc}_{g, 0}^r$
consists of points $[C, x_1, \ldots, x_{(2r+1)(g-1)}]$ such that
$\sum_{i=1}^{(2r+1)(g-1)} x_j\in |K_C^{\otimes r}|$.
\begin{theorem}\label{mrc}
When $n=(2r+1)(g-1)-2i$, the locus $\mathfrak{Mrc}_{g,i}^r$ is a
divisor on $\cM_{g,n}$ and the class of its compactification in
$\mm_{g, n}$ is given by the following formula:
$$\overline{\mathfrak{Mrc}}_{g,i}^r\equiv \frac{1}{g-1}{g-1\choose
i}\Bigl(a\lambda+c\sum_{j=1}^n \psi_j-b_{0}\delta_{0}-\sum_{j, s\geq
0, } b_{j:s}\sum_{|S|=s} \delta_{j:S}\Bigr),$$ where
$$c=rg+g-i-r-1, \ b_{0}=-\frac{1}{g-2}\Bigl({r+1\choose 2}(g-1)(g-2)+i(i+1+2r-rg-g)\Bigr),$$
$$a=-\frac{1}{g-2}\Bigl((g-1)(g-2)(6r^2+6r+1)+i(24r+10i+10-10g-12rg)\Bigr),$$
$$
\ b_{0:s}={s+1\choose 2}(g-1)+s(rg-r)-si, \mbox{ and } b_{j:s}\geq
b_{0: s} \mbox{ for } j\geq 1.$$
\end{theorem}

Using (\ref{logan}) and Theorem \ref{mrc} one obtains the following
table for which $\mm_{g, n}$'s are known to be of general type. In
each case the strategy is to show that $K_{\mm_{g, n}}$ lies in the
cone spanned by $\overline{D}_{g: a_1, \ldots, a_n}$, (pullbacks of)
$\overline{\mathfrak{Mrc}}_{g, i}^r$, boundary divisors and ample
classes:

\begin{theorem}\label{mgn}
For integers $g=4, \ldots, 21$, the moduli space $\mm_{g, n}$ is of
general type for all $n\geq f(g)$ where $f(g)$ is described in the
following table.
\begin{center}
\begin{tabular}{c|cccccccccccccccccc}
$g$ & 4& 5& 6& 7& 8& 9& 10& 11& 12& 13& 14& 15& 16& 17& 18& 19& 20&
21\\
\hline $f(g)$ & 16 & 15&16& 15 &14 & 13& 11& 12& 13& 11& 10& 10& 9&
9&
9& 7& 6& 4\\
\texttt{}\end{tabular}
\end{center}
\end{theorem}

We end this paper with a number of questions related to the global
geometry of $\mm_g$.

\subsection{The hyperbolic nature of $\mm_g$}

Since $\mm_g$ is of general type for $g$ large, through a general
point there can pass no rational or elliptic curve (On the other
hand there are plenty of rational curves in $\mm_g$ for every $g$,
take any pencil on a surface).

\begin{question}
For large $g$, find a suitable lower bound for the invariant
$$\gamma_g:=\mbox{inf}\{g(\Gamma): \Gamma \subset \mm_g \mbox{ is a
curve passing through a general point } [C]\in \mm_g\}.$$ Is it true
that for large $g$ we have that $\gamma_g\geq C\ \mbox{log}(g)$,
where $C$ is a constant independent of $g$? At present we do not
even seem to be able to rule out the (truly preposterous)
possibility that $\gamma_g=2$. Since $\Gamma$ will correspond to a
fibration $f:S\rightarrow \Gamma$ from a surface with fibres being
curves of genus $g$, using the well-known formulas
$$\Gamma \cdot \delta=c_2(S)+4(g-1)(1-g(\Gamma)) \mbox{ and }\Gamma
\cdot \lambda=\chi(\OO_S)+(g-1)(1-g(\Gamma)),$$ the question can
easily be turned into a problem about the irregularity of surfaces.
\end{question}

\begin{question}
For large $g$, compute the invariant
$$n_g:=\mbox{max}\{\mbox{dim}(Z): Z\subset \mm_g, Z\cap \cM_g\neq
\emptyset , \ Z \mbox{ is not of general type} \}.$$ For the moduli
space $\mathcal{A}_g$, Weissauer proved that for $g\geq 13$ every
subvariety of $\mathcal{A}_g$ of codimension $\leq g-13$ is of
general type (cf. \cite{W}). It seems reasonable to expect something
along the same lines  for $\mm_g$.
\end{question}

\begin{question} A famous theorem of Royden implies that $\cM_{g, n}$
admits no non-trivial automorphisms or unramified correspondences
for $2g-2+n\geq 3$ (see e.g. \cite{M} Theorem 6.1 and the references
cited therein). Recall that a non-trivial unramified correspondence
is a pair of distinct finite \'etale morphisms $\alpha:X\rightarrow
\cM_{g, n}, \beta:X\rightarrow \cM_{g, n}$. Precisely, Royden proves
that the group of holomorphic automorphisms of the Teichm\"uller
space $\mathcal{T}_{g, n}$ is isomorphic to the mapping class group.
Using \cite{GKM} Corollary 0.12 it follows that this result can be
extended to $\mm_{g}$ when $g\geq 1$. Any automorphism
$f:\mm_{g}\rightarrow \mm_{g}$ maps the boundary to itself, hence
$f$ induces an automorphism of $\cM_{g}$ and then $f=1_{\mm_{g,
n}}$. One can ask the following questions:
\begin{enumerate}
\item Is there an algebraic proof of Royden's Theorem (in
arbitrary characteristic) using only intersection theory on $\mm_{g,
n}$?
\item Are there any non-trivial \emph{ramified} correspondences of
$\cM_{g, n}$?
\item Are there any non-trivial \emph{birational} automorphisms of
$\cM_g$ for $g\geq 3$? It is known that there exists an integer
$g_0$ such that $\cA_g$ admits no birational automorphisms for any
$g\geq g_0$ (cf. \cite{Fr}).
\item Is it true that for $n\geq 5$ we have that $\mbox{Aut}(\mm_{0,
n})=S_n$? Note that it follows from \cite{KMcK} Theorem 1.3, that
$\mbox{Aut}\bigl(\mm_{0, n}/S_n\bigr)=\{\mathrm{Id}\}$.
\end{enumerate}
\end{question}

\begin{question} We fix an ample line bundle $L\in
\mbox{Pic}(\mm_g)$ (say $L=\kappa_1$). Can one attach a modular
meaning to the \emph{Seshadri constant} $\epsilon_{\mm_g}(L, [C])$,
where $[C]\in \cM_g$ is a general curve?
\end{question}

\end{document}